\documentclass[preprint]{elsarticle1}
\usepackage[english]{babel}
\usepackage{amsmath,amssymb,amsfonts,amsthm,mathrsfs,mathtools}
\usepackage{indentfirst}
\usepackage{enumitem}
\usepackage[autostyle,threshold=0]{csquotes}
\usepackage{natbib}
\usepackage[colorlinks=true]{hyperref}

\usepackage{subcaption}

\newtheorem{proposition}{Proposition}[section]
\newtheorem{remark}{Remark}[section]

\numberwithin{equation}{section}

\journal{}

\begin{document}

\begin{frontmatter}
\title{A negative answer to a conjecture arising in the study \\ of selection-migration models in population genetics\tnoteref{t1}}
\tnotetext[t1]{This work was supported by the auspices of
Grup\-po Na\-zio\-na\-le per l'Anali\-si Ma\-te\-ma\-ti\-ca, la Pro\-ba\-bi\-li\-t\`{a} e le lo\-ro
Appli\-ca\-zio\-ni (GNAMPA) of the Isti\-tu\-to Na\-zio\-na\-le di Al\-ta Ma\-te\-ma\-ti\-ca (INdAM).}
\author{Elisa Sovrano}
\address{Department of Mathematics, Computer Science and Physics, University of Udine, \\via delle Scienze 206, 33100 Udine, Italy }

\begin{abstract}
We deal with the study of the evolution of the allelic frequencies, at a single locus, 
for a population distributed continuously over a bounded habitat. 
We consider evolution which occurs under the joint action of selection and arbitrary migration, 
that is independent of genotype, in absence of mutation and random drift. 
The focus is on a conjecture, that was raised up in literature of population genetics, 
about the possible uniqueness of polymorphic equilibria, which are known as clines, 
under particular circumstances. 
We study the number of these equilibria, making use of topological tools, 
and we give a negative answer to that question by means of two examples. 
Indeed, we provide numerical evidence of multiplicity of positive solutions
for two different Neumann problems satisfying the requests of the conjecture.
\end{abstract}

\end{frontmatter}

\section{Introduction}
Population genetics is a field of biology concerning the genetic structure inside the populations. 
Its main interest is the understanding of evolutionary processes that make the complexity of the Nature so intriguing. 
One of the main causes of the diversity among organisms are the changes in the genetic sequence. 
The genome evolution is influenced by selection, recombination, harmful and beneficial mutations, among others. 
This way, population genetics becomes helpful in order to tackle a very broad class of issues, 
from epidemiology, animal or plant breeding, demography and ecology.

The birth of the ``modern population genetics'' can be traced back at the need to interlace the Darwin's evolution theory with the Mendelian laws of inheritance. 
This has taken place in the 1920s and early 1930s, when Fisher, Haldane and Wright 
have developed mathematical models in order to analyze how the natural selection, along with other factors, 
would modify the genetic composition of a population over time. 
Accordingly, an impressive moment on the history of this field of the genetics is the ``Sixth International Congress of Genetics'' in Ithaca, 
where all the three fathers of the genetical theory of evolution have given a presentation of their pioneering works, see \cite{proceeding-32}.

Mathematical models of population genetics can be described by relative genotypic frequencies or 
relative allelic frequencies, that may depend on both space and time. 
A common assumption is that individuals mate at random in a habitat 
(which can be bounded or not) with respect to the locus under consideration. 
Furthermore, the population is usually considered large enough so that frequencies can be treated as deterministic. 
This way, a probability is associated to the relative frequencies of genotypes/alleles.
The dynamics of gene frequencies are the result of some genetic principles along with several environmental influences, 
such as selection, segregation, migration, mutation, recombination and mating, that lead to different evolutionary processes 
like adaptation and speciation, see \cite{Bu-14}. 

Amongst these influences, by \emph{natural selection} we mean that some genotypes enjoy a survival or reproduction advantage over other ones. 
This way, the genotypic and allelic frequencies change in accord to the proportion of progeny to the next generation of 
the various genotypes which is named \emph{fitness}. 
Thinking to model real-life populations, we have to take into account which is unusual that the selection factor acts alone. 
Since every organism lives in environments that are heterogeneous, 
another considerable factor is the natural subdivision of the population that mate at random only locally. 
Thus, \emph{migration} is often considered as a factor that affects the amount of genetic change. 
There are two different ways, in order to model the dispersion or the migration of organisms: 
one is of discrete type and the other one of a continuous nature. 
If the population size is sufficiently large and the selection is restricted to a single locus with two alleles, 
then deterministic models continuous in time and space lead to mathematical problems 
which involve a single nonlinear partial differential equation of reaction-diffusion type.

In this direction, a seminal paper was given by \cite{Fi-37}. 
In that work, it was studied the frequency of an advantageous gene for a uniformly distributed population in a one-dimensional habitat 
which spreads through an intensity constant selection term. 
Accordingly, a mathematical model of a \emph{cline} was built up as a \emph{non-constant stationary solution} 
of the nonlinear diffusion equation in question. 
The term cline was coined by J. Huxley in \cite{Hu-38}: 
``Some special term seems desirable to direct attention to variation within groups, and I propose the term cline, 
meaning a gradation in measurable characters.''
One of the major causes of cline's occurrence is the migration or the selection which favors an allele in a region of the habitat 
and, a different one in another region. 
The \emph{steepness of a cline} is considered as an indicating character of the level of the geographical variation.
Another contribution comes from J.B. Haldane in \cite{Ha-48}, who has studied the cline's stability by considering 
as a selection term a stepwise function which depends on the space and changes its sign. 

Some meaningful generalizations of these models have been performed, for example, 
in \cite{Fi-50} by introducing a linear spatial dependence in the selection term; in
\cite{Sl-73} by considering a different diffusion term that can model barriers and in 
\cite{Na-75} by taking into account population not necessarily uniformly distributed and 
terms of migration-selection that depend on both space and time.
During the past decade, these mathematical treatments have opened the door to a great amount of works that 
investigated the existence, uniqueness and stability of clines. 
Since a complete list of references of further analysis on clines is out of the scope of this work, 
we limit ourselves to cite some of the earliest contributions in the literature that have inspired the succeeding ones, 
see for instance \cite{Co-74,FiPe-77,Fl-75,Na-76,Na-78,Pe-78}.

Understanding the processes that act in order to have non-constant genetic polymorphisms 
(i.e., loci that occur in more than one allelic form) is an important challenge in population genetics. 
In the present work, we deal with a class of diallelic migration-selection models in continuous space and time 
introduced by W.H. Fleming in \cite{Fl-75} and D. Henry in \cite{He-81}.
We focus on a conjecture stated in \cite{LoNa-02},
that, for such a kind of reaction diffusion equations, guess the uniqueness of a cline (instead of the existence of multiple ones). 
In a one-dimensional setting, we will give a negative answer to that conjecture, 
by providing two examples with multiplicity of non-constant steady states. 
This type of treatment is inspired by the result, about multiplicity of positive solutions for indefinite weight problems 
with Dirichlet boundary conditions, performed in \cite{SoZa-15}.
Although the problem approach has a topological feature, numerical simulations are given in order to support it.

The plan of the paper is the following. 
In Section~\ref{section-2}, we present the class of migration-selection models 
considered and the state of the art which has lead to the formulation of the conjecture of Lou and Nagylaki,
with reference to the genetical and mathematical literature.
In Section~\ref{section-3}, we build up two examples giving a negative answer to this conjecture.
In Section~\ref{section-4}, we conclude with a discussion.

\section{Migration-selection model: the conjecture of Lou and Nagylaki}
\label{section-2}
\def\theequation{2.\arabic{equation}} \makeatother
To ease understanding the conjecture raised up in \cite{LoNa-02}, 
we introduce some notations. 
We also provide an overview of the classical migration-selection model, continuous in space and in time,
of a population in which the genetic diversity occurs in one locus with two alleles, $A_{1}$ and $A_{2}$. 

Let us consider a population continuously distributed in a bounded habitat, say $\Omega$. 
In our context, genetic diversity is the result only of the joint action of dispersal within $\Omega$ 
and selective advantage  for some genotypes, 
so that, no mutation nor genetic drift will be considered. 
This way, the gene frequencies, after random mating, are given by the Hardy-Weinberg relation. 
The genetic structure of the population is measured by the frequencies $p(x,t)$ and 
$q(x,t):=(1-p(x,t))$ at time $t$ and location $x\in\Omega$ of $A_{1}$ and $A_{2}$, respectively. 
Thus, by the assumptions made, the mathematical formulation of this migration-selection model 
leads to the following semilinear parabolic PDE:
\begin{equation}\label{eq-pde}
\frac{\partial p}{\partial t} = \Delta p + \lambda w(x) f(p) \quad \text{ in }\Omega\times]0,\infty[,
\end{equation}
where $\Delta$ denotes the Laplace operator and $\Omega \subseteq {\mathbb R}^N$ is a bounded open connected set, with $N\geq1$,
whose boundary $\partial \Omega$ is $C^2$.
The term $\lambda w(x) f(u)$ models the effect of the natural selection. 
More in detail, the real parameter $\lambda > 0$ plays the role of the ratio of the selection intensity and
the function $w\in L^{\infty}(\Omega)$ represents the local selective advantage (if $w(x) > 0$), 
or disadvantage (if $w(x) < 0$), of the gene at the position $x\in \Omega.$ 
Moreover, following \cite{Fl-75} and \cite{He-81}, the nonlinear term we treat is a general function $f: [0,1]\to {\mathbb R}$ of class $C^2$ satisfying 
\begin{equation*}
f(0) = f(1) =0, \quad f(s) > 0 \;\; \forall\,s\in \,]0,1[\,, \quad f'(0) > 0 > f'(1).
\leqno{(f_{*})}
\end{equation*}
We also impose that there is no-flux of genes into or out of the habitat $\Omega$, 
namely we assume that
\begin{equation}\label{eq-boundary}
\frac{\partial p}{\partial \nu} = 0 \quad \text{ on }\partial\Omega\times]0,\infty[,
\end{equation}
where $\nu$ is the outward unit normal vector on $\partial\Omega$.
Since $p(t,x)$ is a frequency, then we are interested only in positive solutions of \eqref{eq-pde}--\eqref{eq-boundary} such that $0\leq p\leq 1$.

By the analysis developed in \cite{He-81}, 
we know that, if the conditions in $(f_{*})$ hold and $0\leq p(\cdot,0) \leq 1$ in $\Omega$, then $0\leq p(x,t) \leq 1$ 
for all $(x,t)\in\Omega\times]0,\infty[$ and equation \eqref{eq-pde} defines a dynamical system in 
\[
X:=\{p\in H^1(\Omega): 0 \leq p(x) \leq 1, \; \text{a.e. in } \Omega\},
\]
where $H^{1}(\Omega)$ is the standard Sobolev space of integrable functions whose first derivative is also square integrable.
Moreover, the stability of the solutions is determined by the equilibrium solutions in the space $X.$ 
Clearly, a stationary solution of the problem \mbox{\eqref{eq-pde}--\eqref{eq-boundary}} is a function $p(\cdot)$ satisfying $0 \leq p\leq 1$, 
\begin{equation}\label{eq-ode}
-\Delta p = \lambda w(x) f(p) \quad \text{ in }\Omega
\end{equation}
and the Neumann boundary condition 
\begin{equation}\label{eq-boundary2}
\frac{\partial p}{\partial \nu} = 0 \quad \text{ on }\partial\Omega.
\end{equation}

Notice that $p\equiv 0$ and $p\equiv 1$ are constant trivial solutions 
to problem \mbox{\eqref{eq-ode}--\eqref{eq-boundary2}}, that correspond to monomorphic equilibria, 
namely when, in the population, the allele $A_{2}$ or $A_{1}$, respectively, is gone to fixation. 
So, one is interested in finding non-constant stationary solutions or, in other words, polymorphic equilibria. 
Indeed, our main interest is the existence of \emph{clines} for system \eqref{eq-pde}--\eqref{eq-boundary}. 

The maintenance of genetic diversity is examined by seeking for the existence 
of polymorphic stationary solutions/clines, 
that are solutions $p(\cdot)$ to system 
\begin{equation*}
\begin{dcases}
-\Delta p = \lambda w(x) f(p) & \text{ in }\Omega,\\
\frac{\partial p}{\partial \nu} = 0 & \text{ on }\partial\Omega,
\end{dcases}
\eqno{(\mathscr{N}_{\lambda})}
\end{equation*}
with $0 < p(x) < 1$ for all $x\in\overline{\Omega}$.

In this respect, the assumption $f(s) > 0$ for every $s > 0$ implies that a necessary condition 
for positive solutions of problem $(\mathscr{N}_{\lambda})$ is that the function $w$ changes its sign. 
In fact, by integrating \eqref{eq-ode} over $\Omega$, we obtain
\[
0=\int_{\Omega}\Delta p +\lambda w(x) f(p)\, dx=
\int_{\partial\Omega}\frac{\partial p}{\partial \nu} \,dx+\lambda\int_{\Omega} w(x) f(p)\,dx =
\lambda\int_{\Omega} w(x) f(p)\,dx. 
\]
Notice that we can see the function $w$ in $(\mathscr{N}_{\lambda})$ 
as a weight term which attains both positive and negative values, 
so that such a kind of system is usually known as \emph{problem with indefinite weight}. 

It is a well-known fact that the existence of positive solutions of $(\mathscr{N}_{\lambda})$ 
depends on the sign of 
\begin{equation}\label{eq-w}
\bar{w}:=\int_{\Omega}w(x)\,dx.
\end{equation}
Indeed, for the linear eigenvalue problem $-\Delta p(x)=\lambda w(x)p(x)$, 
under Neumann boundary condition on $\Omega$, the following facts hold:
if $\bar{w}<0$, then there exists a unique positive eigenvalue having an associated eigenfunction which does not change sign; 
on the contrary, if $\bar{w}\geq0$ such an eigenvalue does not exist and $0$ is the only non-negative eigenvalue 
for which the corresponding eigenfunction does not vanish, see \cite[Theorem~3.13]{BrLi-80}.

Furthermore, under the additional assumption of concavity for the nonlinearity:
\begin{equation}\label{eq-concave}
f''(s)\leq0, \quad \forall s>0,
\end{equation}
it follows that, if $\bar{w}<0,$
then there exists $\lambda_{0}>0$ such that for each $\lambda>\lambda_{0}$ 
problem \mbox{\eqref{eq-pde}--\eqref{eq-boundary}} has at most one nonconstant stationary positive solution 
(i.e., cline) which is asymptotically stable, see \cite[Theorem~10.1.6]{He-81}.

After these works a great deal of contributions appeared in order 
to complement these results of existence and uniqueness on population genetics, see for instance \cite{BaPoTe-88,BrHe-90,BLT-89}; 
or to consider also unbounded habitats as done in \citep{FiPe-81}; 
or even to treat more general uniformly elliptic operators, as in \cite{Se-83,SeHe-82}.
Taking into account these works, in \cite{LoNa-02} the migration-selection model 
with an isotropic dispersion, that is identified with the Laplacian operator,
was generalized to an arbitrary migration, which involves a strongly uniformly elliptic differential operator of second order 
(see also \cite{Na-89,Na-96} for the derivation of this model as a continuous approximation of the discrete one).

By modeling single locus diallelic populations, there is an interesting family of nonlinearities 
which satisfies the conditions in $(f_{*})$ and 
allows to consider different phenotypes of alleles, $A_{1}$ and $A_{2}$. 
This family can be obtained by considering
the map $f_{k}:\mathbb{R}^{+}\to\mathbb{R}^{+}$ such that
\begin{equation}\label{eq-linearity}
f_{k}(s):= s(1-s)(1+k-2ks),
\end{equation}
where $-1\leq k\leq1$ represents the degree of dominance of the alleles independently of the space variable, see \cite{Na-75}.
In this special case, if $k=0$ then the model does not present any kind of dominance, 
instead, if $k=1$ or $k=-1$ then the allelic dominance is relative to $A_{1}$, in the first case, and to $A_{2}$ in the second one 
(the last is also equivalent to said that $A_{1}$ is recessive).
In view of this, we can make mainly the following two observations.

\begin{remark}\label{rem-1}
In the case of no dominance, i.e. $k=0$, from \eqref{eq-linearity} we have $f_{0}(s)=s(1-s)$ which is a concave function. 
Therefore, we can enter in the settings considered by \emph{\cite{He-81}.}
So if $w(x)>0$ on a set of positive measure in $\Omega$ and $\bar{w}<0$, 
then for $\lambda$ sufficiently large there exists a unique positive non trivial equilibrium of the equation
$\partial p/ \partial t = \Delta u + \lambda w(x) p(1-p)$ for every $(x,t)\in\Omega\times]0,\infty[$ under 
the boundary condition \eqref{eq-boundary}.
\end{remark}

\begin{remark}\label{rem-2}
In the case of completely dominance of allele $A_{2}$, i.e. $k=-1$, from \eqref{eq-linearity}
we have $f_{-1}(s)=2s^{2}(1-s)$ which is not a concave function. 
Thanks to the results in \emph{\cite{LNs-10}}, if $w(x)>0$ on a set of positive measure in $\Omega$ and $\bar{w}<0$, 
then for $\lambda$ sufficiently large there exist at least two positive non trivial equilibrium of the equation
$\partial p/ \partial t = \Delta p + \lambda w(x) 2 p^{2}(1-p)$ for every $(x,t)\in\Omega\times]0,\infty[$ under 
the boundary condition \eqref{eq-boundary}.
\end{remark}

We observe that the map $s\mapsto f_{0}(s)/s$ is \emph{strictly decreasing} with $f_{0}(s)$ \emph{concave}. 
On the contrary, the map $s\mapsto f_{-1}(s)/s$ is \emph{not strictly decreasing} with $f_{-1}(s)$ \emph{not concave}. 
Thus, from Remark~\ref{rem-1} and Remark~\ref{rem-2}, it arises a natural question which involves the possibility 
to weaken the concavity assumption \eqref{eq-concave} further to the monotonicity of the map $s\mapsto f(s)/s$, in order to get 
uniqueness results of nontrivial equilibria for problem \mbox{\eqref{eq-pde}--\eqref{eq-boundary}}. 

This is still an open question, firstly appeared in \citep[Conjecture~5.1]{LoNa-02}, known as the ``conjecture of Lou and Nagylaki''.

\medskip

\noindent \textit{Conjecture} ``\emph{Suppose that $w(x)>0$ on a set of positive measure in $\Omega$ and such that $\bar{w}=\int_{\Omega} w \,dx<0$. 
If the map $s\mapsto f(s)/s$ is monotone decreasing in $]0,1[$, 
then \mbox{\eqref{eq-pde}--\eqref{eq-boundary}} has at most one nontrivial equilibrium $p(t,x)$ with $0<p(0,x)<1$ for every $x\in\bar{\Omega}$,
which, if it exists, is globally asymptotically stable.}'' \citep[c.f.][p.~4364]{LNN-13}.

\medskip

\noindent The study of existence, uniqueness and multiplicity of positive solutions for nonlinear indefinite weight problems 
is a very active area of research, starting from the Seventies, 
and several types of boundary conditions along with a wide variety of nonlinear functions, classified according to growth conditions, 
were taken into account.
Several authors have addressed this topic, see \cite{AlTa-93,AmLG-98,BCN-94,BCN-95,BrOs-86,HeKa-80}, 
just to recall the first main papers dedicated.
Instead, the recent literature about multiplicity results for positive solutions of indefinite weight problems with 
Dirichlet or Neumann boundary conditions is really very rich. 
In order to cover most of the results achieved with different techniques so far, 
we give reference of the following bibliography \cite{BGH-05,Bo-11,BGZ-16,BoGa-16,FeZa-15,GHZ-03,GiGo-09,GoLo-00,ObOm-06}. 

Nevertheless, as far as we known, there is no answer about the conjecture of Lou and Nagylaki.
It is interesting to notice that the study of the concavity of $f(s)$ versus the monotonicity of $f(s)/s$ has significance also 
in the investigation on the uniqueness of positive solutions for a particular class of indefinite weight problems with Dirichlet boundary conditions.
More in detail, these problems involve positive nonlinearities which have linear growth at zero and sublinear growth at infinity,
namely
\begin{equation*}
\begin{dcases}
-\Delta p = \lambda w(x) g(p) & \text{ in }\Omega,\\
p = 0 & \text{ on }\partial\Omega,
\end{dcases}
\eqno{(\mathscr{D}_{\lambda})}
\end{equation*}
where $g:\mathbb{R}^{+}\to\mathbb{R}^{+}$ is a continuous function satisfying 
\begin{equation*}
g(0)=0\,,\quad g(s) > 0 \;\; \forall\,s>0\,, \quad \lim_{s\to0^{+}}\frac{g(s)}{s}>0\,,\quad \lim_{s\to+\infty}\frac{g(s)}{s}=0.
\leqno{(g_{*})}
\end{equation*}

The state of the art on this topic refers mainly on two papers.
From the results achieved in \cite{BrOs-86}, it follows that, 
if $(g_{*})$ holds and, moreover, the map $s\to g(s)/s$ is strictly decreasing, 
then there exist at most one positive solution of $(\mathscr{D}_{\lambda})$ only 
if the weight function $w(x) > 0$ for a.e. $x\in\Omega.$
On the other hand, from \cite{BrHe-90}, if the conditions in $(g_{*})$ are satisfied for a smooth concave nonlinear term $g$  
and the weight $w$ is a smooth and changing sign function, then there exists at most one positive solutions of $(\mathscr{D}_{\lambda})$.
Therefore, if the weight function is positive, then the hypothesis of Brezis-Oswald, concerning the monotonicity of $g(s)/s$, 
is more general than the requirement of Brown-Hess about the concavity of $g(s)$.

At this point one could query whether something similar to the conjecture of Lou and Nagylaki happens 
also for this family of Dirichlet problems.
This was done in \citep[ch.~5]{SoZa-15}, where it was shown that the monotonicity of the map $s\mapsto g(s)/s$ is not enough 
to guarantee the uniqueness of positive solutions for problems as in $(\mathscr{D}_{\lambda})$ with an indefinite weight. 
Through numerical evidence, more than one positive solution has been detected for an exemplary 
two-point boundary value problem $(\mathscr{D}_{\lambda})$.

\section{Multiplicity of clines: the conjecture has negative answer}
\label{section-3}
\def\theequation{3.\arabic{equation}} \makeatother
In this section we look at the framework of the conjecture of Lou and Nagylaki. 
So, from now on we tacitly consider a nonlinear function $f: [0,1]\to {\mathbb R}$ of class $C^2$ 
which is \emph{not concave}, satisfies $(f_{*})$ and is such that the map \emph{$s\mapsto f(s)/s$ is strictly decreasing}. 

We concentrate on the one-dimensional case $N=1$ 
and we take as a habitat an open interval $\Omega:=]\omega_{1},\omega_{2}[$ with $\omega_{1},\omega_{2}\in\mathbb{R}$ 
such that $\omega_{1}<0<\omega_{2}$.
This type of habitats, confined to one-dimensional spaces, have an intrinsic interest in modeling phenomena which occur, 
for example, in neighborhoods of rivers, sea shores or hills, see \cite{Na-78}.
As in  \cite{Na-75,Na-78}, we assume that the weight term $w$ is step-wise.
Hence, let us consider the following class of indefinite weight functions 
\begin{equation}\label{eq-weight}
w(x):=
\begin{cases}
-\alpha	& x\in[\omega_{1},0[,\\
1	& x\in]0,\omega_{2}],
\end{cases}
\end{equation}
such that 
\[\bar{w} = -\omega_{1}\alpha +\omega_{2}<0,\] 
with $\bar{w}$ defined as in \eqref{eq-w}.
In these settings, the indefinite Neumann problem $(\mathscr{N}_{\lambda})$ reads as follows
\begin{equation}\label{eq-3.1}
\begin{cases}
p'' + \lambda w(x) f(p) = 0,\\
p'(\omega_{1})=0=p'(\omega_{2}),
\end{cases}
\end{equation}
with $0 < p(x) < 1$ for all $x\in[\omega_{1},\omega_{2}].$

Inspired by the results in \citep[ch.~5]{SoZa-15}, we will consider two particular functions $f$ 
in order to provide a negative reply to the conjecture under examination.
In both cases, we are going to use a topological argument, that is called \emph{shooting method}, 
and, with the aid of some numerical computations, we give evidence
of multiplicity of positive solutions for the corresponding problems in \eqref{eq-3.1}.
The shooting method relies on the study of the deformation of planar continua under the action of the vector field associated to 
the second order scalar differential equation in \eqref{eq-3.1}, whose formulation, in the phase-plane $(u,v)$,
is equivalent to the first order planar system
\begin{equation}\label{eq-3.2}
\begin{cases}
u '=v, 	\\
v '= -\lambda w(x) f(u).
\end{cases}
\end{equation}
Solutions $p(\cdot)$ of problem \eqref{eq-3.1} we are looking for are also solutions $(u(\cdot),v(\cdot))$ of system \eqref{eq-3.2},
such that $v(\omega_{1})=0=v(\omega_{2})$.

We set the interval $[0,1]$ contained in the $u-$axis as follows
\[\mathcal{L}_{\{v=0\}}:=\{(u,v)\in\mathbb{R}^{2}: 0\leq u\leq 1,\, v=0\}.\]
This way, as a real parameter $r$ ranges between $0$ and $1$, we are interested in the solution, 
$(u(\cdot\, ;\omega_{1},(r,0)),v(\cdot\, ;\omega_{1},(r,0)))$, of the Cauchy problem
with initial conditions
\begin{equation}\label{eq-ci}
\begin{cases}
u(\omega_{1})=r,\\
v(\omega_{1})=0,
\end{cases}
\end{equation}
such that $(u(\omega_{2} ;\omega_{1},(r,0)),v(\omega_{2} ;\omega_{1},(r,0)))\in \mathcal{L}_{\{v=0\}}$.
Hence, let us consider the planar continuum $\Gamma$ obtained by shooting $\mathcal{L}_{\{v=0\}}$ forward
from $\omega_{1}$ to $\omega_{2}$, namely
\[
\Gamma:=\{(u(\omega_{2};r),v(\omega_{2};r))\in\mathbb{R}^{2} : r\in[0,1]	\}.
\]
We define the set of the intersection points between this continuum and the segment $[0,1]$ contained in the $u-$axis, as
\[\mathcal{S}:=	\Gamma\cap\mathcal{L}_{\{v=0\}}.
\]
Then, there exists an injection form the set of the solutions $p(\cdot)$ of \eqref{eq-3.1} such that $0 < p(x) < 1$ for all $x\in[\omega_{1},\omega_{2}]$
and the set $\mathcal{S}\setminus \left( \{(0,0)\}\cup  \{(1,0)\} \right)$.

More formally, we denote 
by $\zeta(\cdot\,; \omega_{0},z_{0})=(u(\cdot\,; \omega_{0},z_{0}),v(\cdot\,; \omega_{0},z_{0}))$ 
the solution of \eqref{eq-3.2} with $\omega_{0}\in[\omega_{1},\omega_{2}]$ and initial condition
$\zeta(\omega_{0}\,; \omega_{0},z_{0})=z_{0}=(u_{0},v_{0})\in\mathbb{R}^{2}$.
The uniqueness of the solutions for the associated initial value problems guarantee that 
the Poincar\'e map associated to system \eqref{eq-3.2} is well defined.
Recall that, for any $\tau_{1},\tau_{2}\in[\omega_{1},\omega_{2}]$, the Poincar\'e map for system \eqref{eq-3.2}, 
denoted by $\Phi_{\tau_{1}}^{\tau_{2}}$,
is the planar map which at any point $z_{0}=(u_{0},v_{0})\in\mathbb{R}^{2}$ associates the point $(u(\tau_{2}),v(\tau_{2}))$
where $(u(\cdot),v(\cdot))$ is the solution of \eqref{eq-3.2} with $(u(\tau_{1}),v(\tau_{1}))=z_{0}$.
Notice that $\Phi^{\tau_{2}}_{\tau_{1}}$ is a global diffeomorphism of the plane onto itself.

Under these notations, the recipe of the shooting method is the following.
A solution $p(\cdot)$ of \eqref{eq-3.1} such that $0<p(x)<1$ for all $x\in[\omega_{1},\omega_{2}]$
is identified by a point $(c,0)\in\mathcal{L}_{\{v=0\}}$ whose image through the action of the Poincar\'e map,
say $C:=\Phi_{\omega_{1}}^{\omega_{2}}((c,0))\in\Gamma,$ belongs to $\mathcal{L}_{\{v=0\}}.$
This way, the solution $p(\cdot)$ of the Neumann problem with $p(\omega_{1})=c$ is found 
looking at the first component of the map 
\[x\mapsto\Phi_{\omega_{1}}^{x}((c,0))=(u(x),v(x)),\]
since, by construction, $p'(\omega_{1})=v(\omega_{1})=0$ and $p'(\omega_{2})=v(\omega_{2})=0$.
This means that the set $\mathcal{S}$ is made by points such that, each of them determines univocally 
an initial condition, of the form \eqref{eq-ci}, for which the solution $(u(\cdot),v(\cdot))$ of the Cauchy problem 
associated to \eqref{eq-3.2} verifies $v(\omega_{1})=0=v(\omega_{2}).$

The study of the uniqueness of the clines is based on the study, in the phase plane $(u,v)$, of the qualitative properties 
of the shape of the continuum $\Gamma$ which is the image of $\mathcal{L}_{\{v=0\}}$ 
under the action of the Poincar\'e map $\Phi_{\omega_{1}}^{\omega_{2}}$. 
More in detail, we are interested in find real values $c\in]0,1[$ such that 
\[\Phi^{\omega_{2}}_{\omega_{1}}((c,0))\in\Phi_{\omega_{1}}^{\omega_{2}}\left(\mathcal{L}_{\{v=0\}}\right)\cap\mathcal{L}_{\{v=0\}}.
\]
Indeed, our aim is looking for values $c\in]0,1[$ such that the point $C=\Phi^{\omega_{2}}_{\omega_{1}}((c,0))$ belongs to $\mathcal{S}$.

So, let us show now that there exist Neumann problems as in \eqref{eq-3.2} that admit more than one positive solution.
Namely, there exist more than one polymorphic stationary solution for the equation:
\begin{equation}\label{eq-rd}
\frac{\partial p}{\partial t}=p''+\lambda w(x)f(p).
\end{equation}
Roughly speaking, if $\Gamma$ crosses the $u-$axis more than one time, out of the points $(0,0)$ and $(1,0)$,
then $\#(\mathcal{S}\setminus \left( \{(0,0)\}\cup  \{(1,0)\} \right))>1$ and so, we expect a result of non-uniqueness of clines
for equation \eqref{eq-rd}.

\subsection{First example}
Taking into account the definition of the functions in \eqref{eq-linearity}, given a real parameter $h>0$, let us consider the family of maps 
$\hat{f}_{h}:[0,1]\to\mathbb{R}$ of class $C^{2}$ such that 
\[\hat{f}_{h}(s):=s(1-s)(1-h s+h s^{2}).\] 
By definition $\hat{f}(0)=0=\hat{f}(1).$ Moreover, to have $s\mapsto \hat{f}_{h}(s)/s$ monotone decreasing in $]0,1[$ it is sufficient to assume $0<h\leq 3$. 
If the parameter $h$ ranges in $]0,3]$, then it is straightforward to check that $\hat{f}_{h}$ is not concave and 
$\hat{f}_{h}(s)>0$ for every $s\in]0,1[.$

Let us fix $h=3.$ Then, in this case, all the conditions in $(f_{*})$ are verified and 
$\hat{f}_{3}(s)=s(1-s)(1-3s+3s^{2})$. 
As a consequence, we point out the following result of multiplicity.

\begin{proposition}\label{prop-1}
Let $f:[0,1]\to\mathbb{R}$ be such that 
\begin{equation}\label{eq-f1}
f(s):=s(1-s)(1-3s+3s^{2}).
\end{equation}
Assume $w:[\omega_{1},\omega_{2}]\to \mathbb{R}$ be defined as in \eqref{eq-weight} 
with $\alpha=1$, $\omega_{1}=-0.21$ and $\omega_{2}=0.2$. 
Then, for $\lambda=45$ the problem \eqref{eq-3.1} 
has at least $3$ solutions such that $0 < p(x) < 1$ for all $x\in[\omega_{1},\omega_{2}]$.
\end{proposition}

Notice that $\bar{w}=-0.01<0$, so we are in the hypotheses of the conjecture.
Now we follow the scheme of the shooting method, in order to detect three polymorphic stationary solutions
for the equation \eqref{eq-rd}. 
This approach, with the help of numerical estimates, will enable us to prove Proposition~\ref{prop-1}.

\begin{figure}[htb]
	\centering
	\includegraphics[width=1\textwidth]{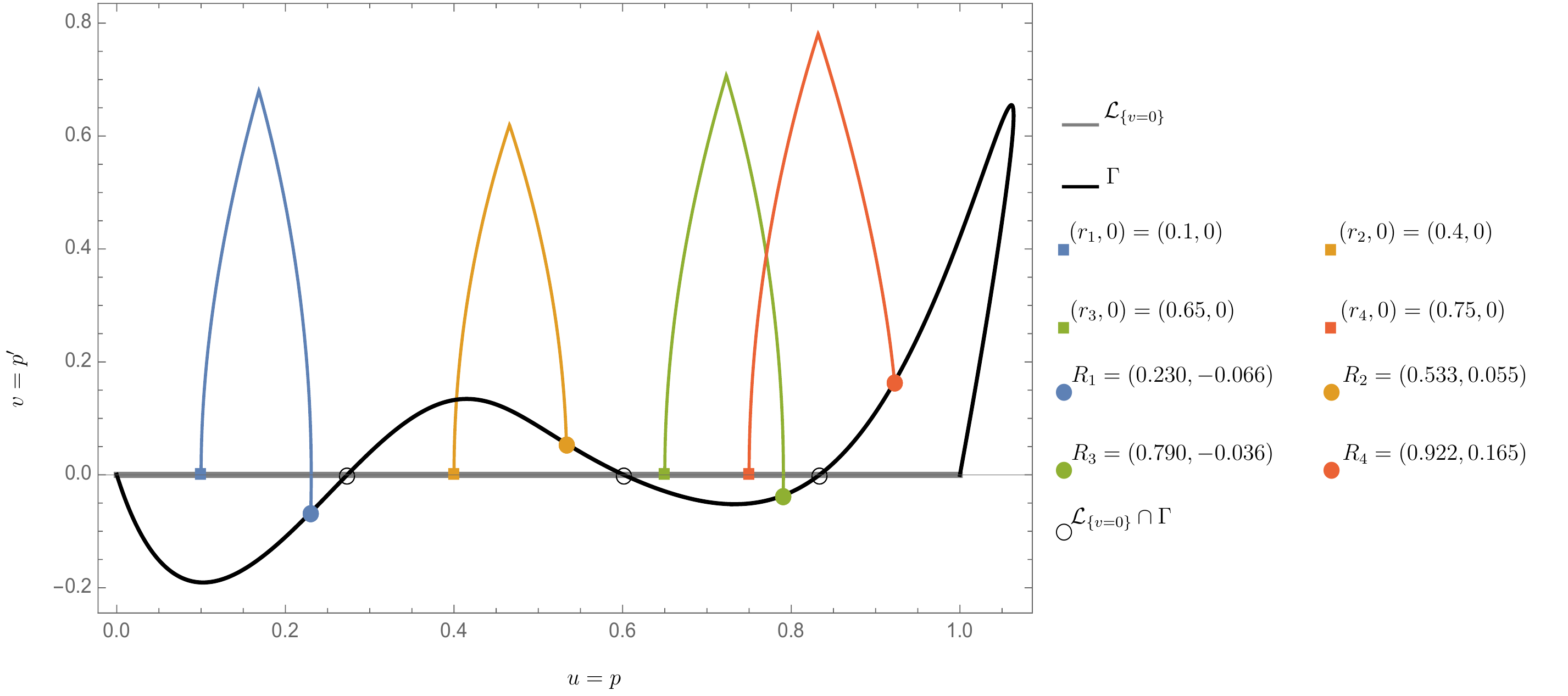}  
\caption{In the phase-plane $(u,v)$: intersection points between $\mathcal{L}_{\{v=0\}}$ and 
$\Gamma=\Phi^{\omega_{2}}_{\omega_{1}}(\mathcal{L}_{\{v=0\}})$; 
solutions of the Cauchy problem with initial conditions given by $(u(\omega_{1}),v(\omega_{1}))=(r_{i},0)$ and 
numerical approximation of the values $R_{i}=(u(\omega_{2}\, ;\omega_{1},(r_{i},0)),v(\omega_{2}\, ;\omega_{1},(r_{i},0)))$ with $i=1,\dots,4.$
The problem's setting $p''+\lambda w(x)f(p)=0$ is defined as in Proposition~\ref{prop-1}.}
\label{fig-1}
\end{figure}
In the phase-plane $(u,v)$, Figure~\ref{fig-1} shows the existence of at least four points 
$(r_{i},0)\in\mathcal{L}_{\{v=0\}}$ with $i=1,\dots,4$ such that, 
by defining their images through the Poincar\'e map $\Phi^{\omega_{2}}_{\omega_{1}}$ as
$R_{i}:=({R_{i}}^{u},{R_{i}}^{v})=\Phi^{\omega_{2}}_{\omega_{1}}((r_{i},0))\in\Gamma$
for every $i\in\{1,\dots,4\}$,
the following conditions  
\[
{R_{i}}^{v}<0	\text{ for }i=1,3,\quad {R_{i}}^{v}>0	\text{ for }i=2,4,
\]
are satisfied. 
This is done, for example, with the choice of the values $r_{1}=0.1$, $r_{2}=0.4$, $r_{3}=0.65$ and $r_{4}=0.75$.
The solutions of the Cauchy problems associated to system \eqref{eq-3.2}, with initial conditions $(r_{i},0)$ for $i=1,\dots,4$, 
assume at $x=\omega_{2}$ the values $R_{1}=(0.230, -0.066)$, $R_{2}= (0.922, 0.165)$, $R_{3}=(0.790,  0.036)$ 
and $R_{4}=(0.533, 0.055)$, truncated at the third significant digit. 
Therefore, we have ${R_{1}}^{v}<0<{R_{2}}^{v}$, ${R_{2}}^{v}>0>{R_{3}}^{v}$ and ${R_{3}}^{v}<0<{R_{4}}^{v}$.
Then, by a continuity argument (that means an application of the Mean Value Theorem), 
there exist at least three real values $c_{1},c_{2}$ and $c_{3}$ such that
\begin{equation}\label{eq-cond1}
r_{j}<c_{j}<r_{j+1}\quad \text{and}\quad
C_{j}:=\Phi^{\omega_{2}}_{\omega_{1}}((c_{j},0))\in\mathcal{S}\setminus \left( \{(0,0)\}\cup  \{(1,0)\} \right),
\end{equation}
for every $j\in\{1,\dots,i-1\}$. So, let us see how to find such values.

\begin{figure}[htb]
	\centering
	\includegraphics[width=1\textwidth]{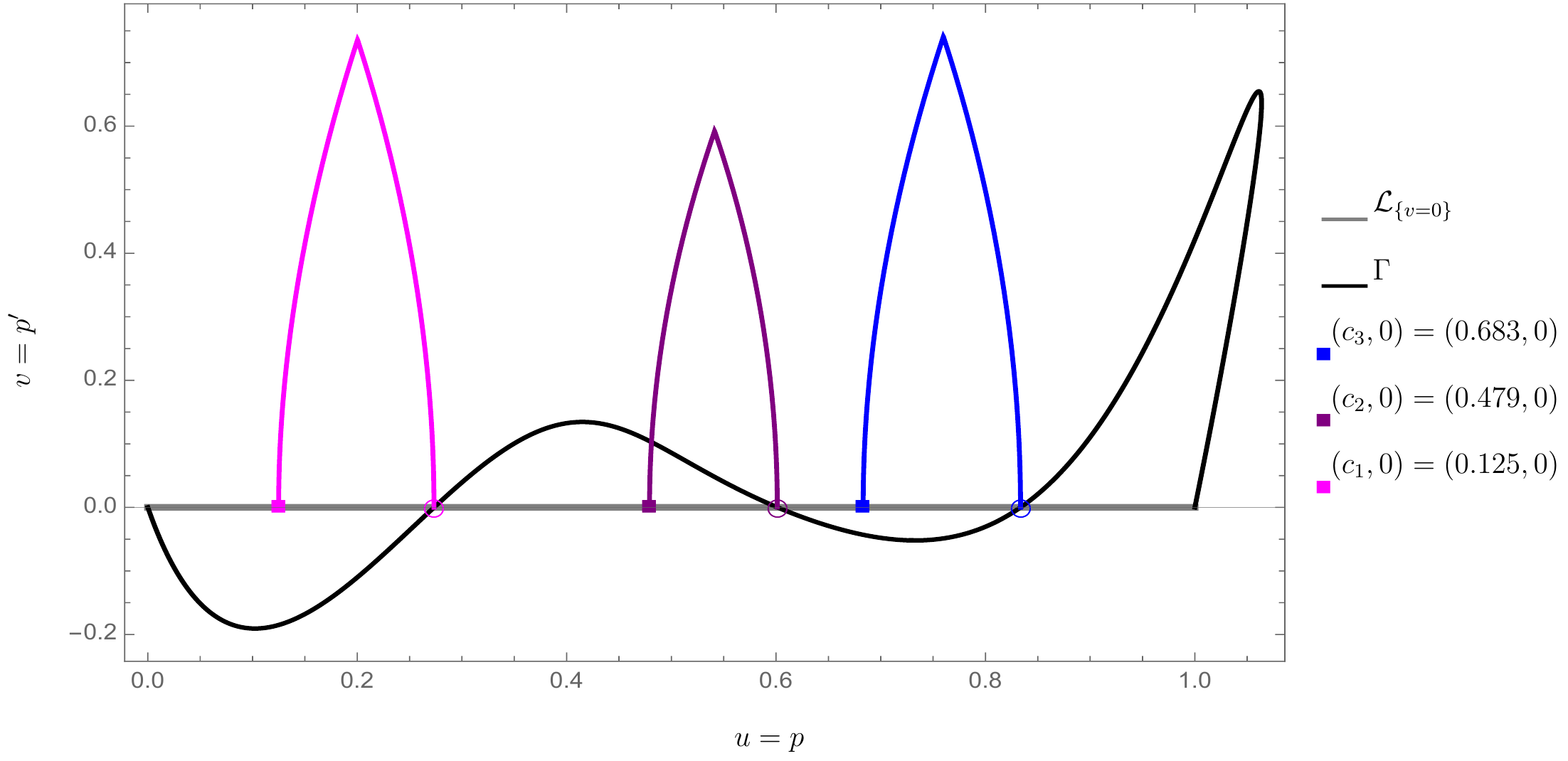}  
\caption{In the phase-plane $(u,v)$: solutions of the Cauchy problem associated to system \eqref{eq-3.2}
with initial conditions given by $(u(\omega_{1}),v(\omega_{1}))=(c_{j},0)$ with $j=1,2,3$. 
The problem's setting $p''+\lambda w(x)f(p)=0$ is defined as in Proposition~\ref{prop-1}.}
\label{fig-2}
\end{figure}

The curve $\Gamma$ is the result of the integration of several system of differential equation \eqref{eq-3.2},
with initial conditions taken within a uniform discretization of the interval $[0,1]$, 
followed by the interpolation of the approximated values of each solution $\zeta(x\,;\omega_{0},z_{0})$ at $x=\omega_{2}$.
Hence, $\Gamma$ represents the approximation of the image of the interval $[0,1]$ 
under the action of the Poincar\'e map $\Phi^{\omega_{2}}_{\omega_{1}}$.

As the Figure~\ref{fig-1} suggests, the projection of $\Gamma$ on its first component 
is not necessarily contained in the interval $[0,1]$, which includes the only values of biological pertinence. 
Nonetheless, this does not avoid the existence of solutions $p(\cdot)$ of the problem \eqref{eq-3.1} such that $0 <p(x)<1$ 
for all $x\in[\omega_{1},\omega_{2}]$.
This way, by means of a fine discretization of $\mathcal{L}_{\{v=0\}}$, 
we have found the approximate values of the intersection points $C_{j}\in 	\Gamma\cap\mathcal{L}_{\{v=0\}}$, 
with $j=1,2,3$. In this case they are:
$C_{1}=(0.273,0)$, $C_{2}=(0.601,0)$ and $C_{3}=(0.833,0)$, truncated at the third significant digit (see Figure~\ref{fig-1}).
The intersection points between $\mathcal{L}_{\{v=0\}}$ and its image $\Gamma$ 
through the Poincar\'e map $\Phi^{\omega_{2}}_{\omega_{1}}$, namely $C_{j}$ with $j=1,2,3$, are in agreement with the previous predictions.

At last, we computed the values $c_{1}=0.125$, $c_{2}=0.479$ and $c_{3}=0.683$, which verify the required conditions \eqref{eq-cond1}. 
For $j=1,2,3$, in Figure~\ref{fig-2} are represented the trajectories of the solutions of the initial value problem
\begin{equation}\label{eq-ivp}
\begin{cases}
p''+\lambda w(x)f(p)=0,\\
p(\omega_{1})=c_{j},\\
p'(\omega_{1})=0,
\end{cases}
\end{equation}
that, by construction, satisfy $p'(\omega_{2})=0$. 
We observe also that the values of each solution $p(\cdot)$ of the three different
initial value problems ranges in $]0,1[$ as desired.

\begin{figure}[htb]
	\centering
	\includegraphics[width=1\textwidth]{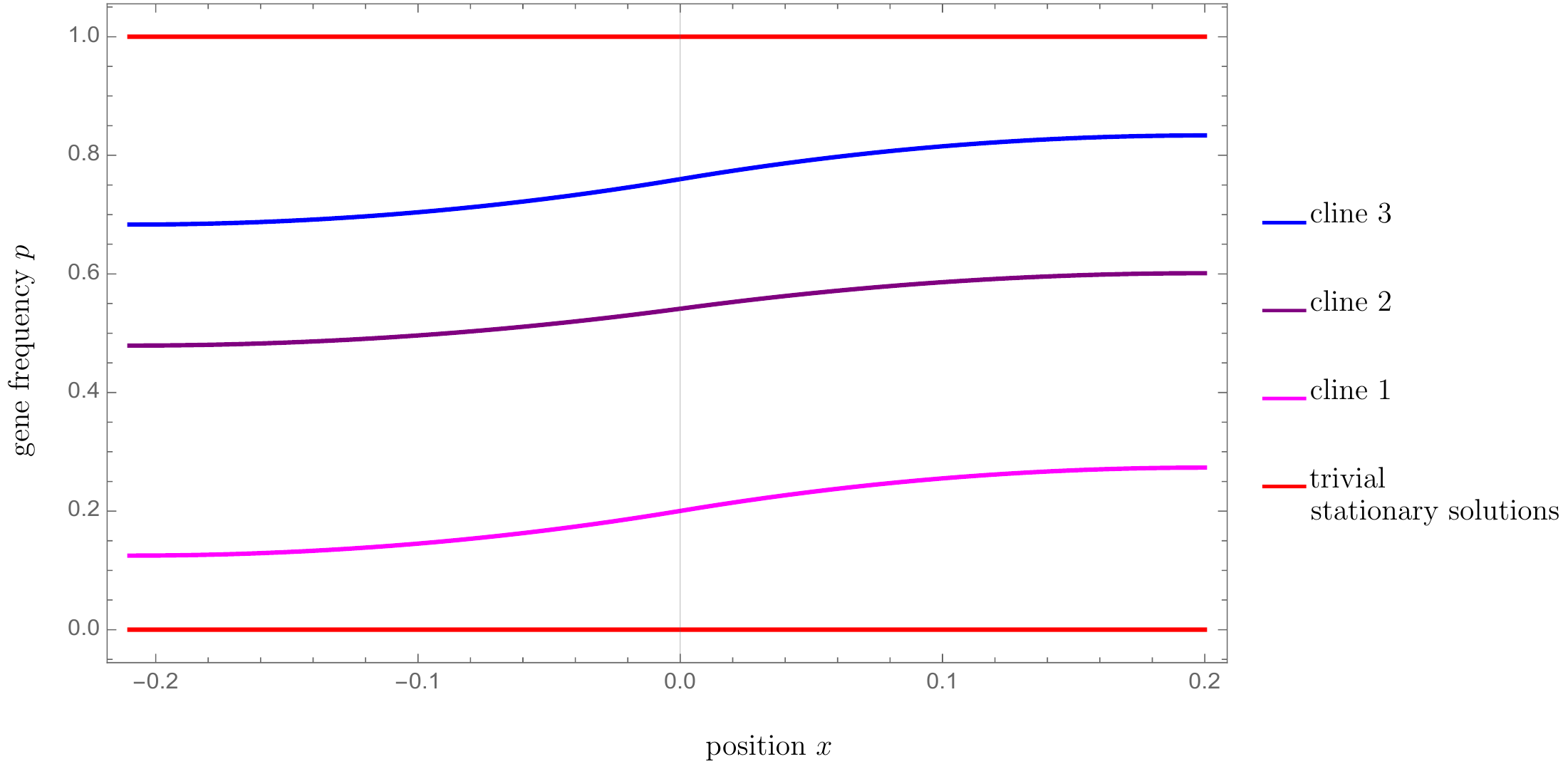}  
\caption{Polymorphic stationary solutions (clines) and trivial stationary solutions ($u\equiv0$ and $u\equiv1$) for equation \eqref{eq-rd}, 
found as positive solutions of the Neumann problem $p''+\lambda w(x)f(p)=0$ satisfying the framework of Proposition~\ref{prop-1}.}
\label{fig-3}
\end{figure}

Once found the values $c_{j}$ with $j=1,2,3$, a numerically result of multiplicity of clines is achieved.
Indeed, in Figure~\ref{fig-3}, we display the approximation of the three nontrivial stationary solution $p(\cdot)$
of equation \eqref{eq-rd} that are identified by the points
$C_{j}\in(\mathcal{S}\setminus \left( \{(0,0)\}\cup  \{(1,0)\} \right))$,
with $j=1,2,3.$

\subsection{Second example}
We refer now to the application given in \citep[ch.~5]{SoZa-15} and we adapt it to our purposes. 
So we consider, the nonlinear term $\tilde{f}:\mathbb{R}^{+}\to\mathbb{R}^{+}$ defined by
\begin{equation*}
\tilde{f}(s):=\left(10 s e^{-25s^{2}}+\frac{s}{|s|+1}\right).
\end{equation*}
It is straightforward to check that $\tilde{f}$ is not concave and the map $s\to \tilde{f}(s)/s$ is strictly decreasing. 
Moreover, $\tilde{f}(0)=0$ and $\tilde{f}(s)>0$ for every $s>0$, but $\tilde{f}$ does not take value zero in $s=1$, 
since $\tilde{f}(1)=10 e^{-25}+1\not=0.$
To satisfy all the conditions in $(f_{*})$, it is sufficient to multiply $\tilde{f}$ by the term $\arctan(m(1-x))$ with $m>0$.
This way, the following result holds.

\begin{proposition}\label{prop-2}
Let $f:[0,1]\to\mathbb{R}$ be such that 
\begin{equation}\label{eq-f2}
f(s):=\left(10 s e^{-25s^{2}}+\frac{s}{|s|+1}\right) \arctan(10-10s).
\end{equation}
Assume $w:[\omega_{1},\omega_{2}]\to \mathbb{R}$ be defined as in \eqref{eq-weight} 
with $\alpha=2.4$, $\omega_{1}=-0.255$ and $\omega_{2}=0.6$. 
Then, for $\lambda=3$ the problem \eqref{eq-3.1} 
has at least $3$ solutions such that $0 < p(x) < 1$ for all $x\in[\omega_{1},\omega_{2}]$.
\end{proposition}

Notice that, under the assumptions of Proposition~\ref{prop-2}, the hypotheses of the conjecture are now all satisfied since $\bar{w}=-0.012<0$. 
To prove the existence of at least three clines for the equation \eqref{eq-rd},
we exploit again the shooting method.

So, our main interest is in finding real values $r_{i}\in]0,1[$ with $i\in\mathbb{N}$ such that,
given $R_{i}:=({R_{i}}^{u},{R_{i}}^{v})=\Phi^{\omega_{2}}_{\omega_{1}}((r_{i},0))$, it follows
\[
{R_{i}}^{v}<0	\text{ for }i=2\ell +1,\quad
{R_{i}}^{v}>0	\text{ for }i=2\ell, \text{ with } \ell\in\mathbb{N}.
\]

\begin{figure}[htb]
	\centering
	\includegraphics[width=1\textwidth]{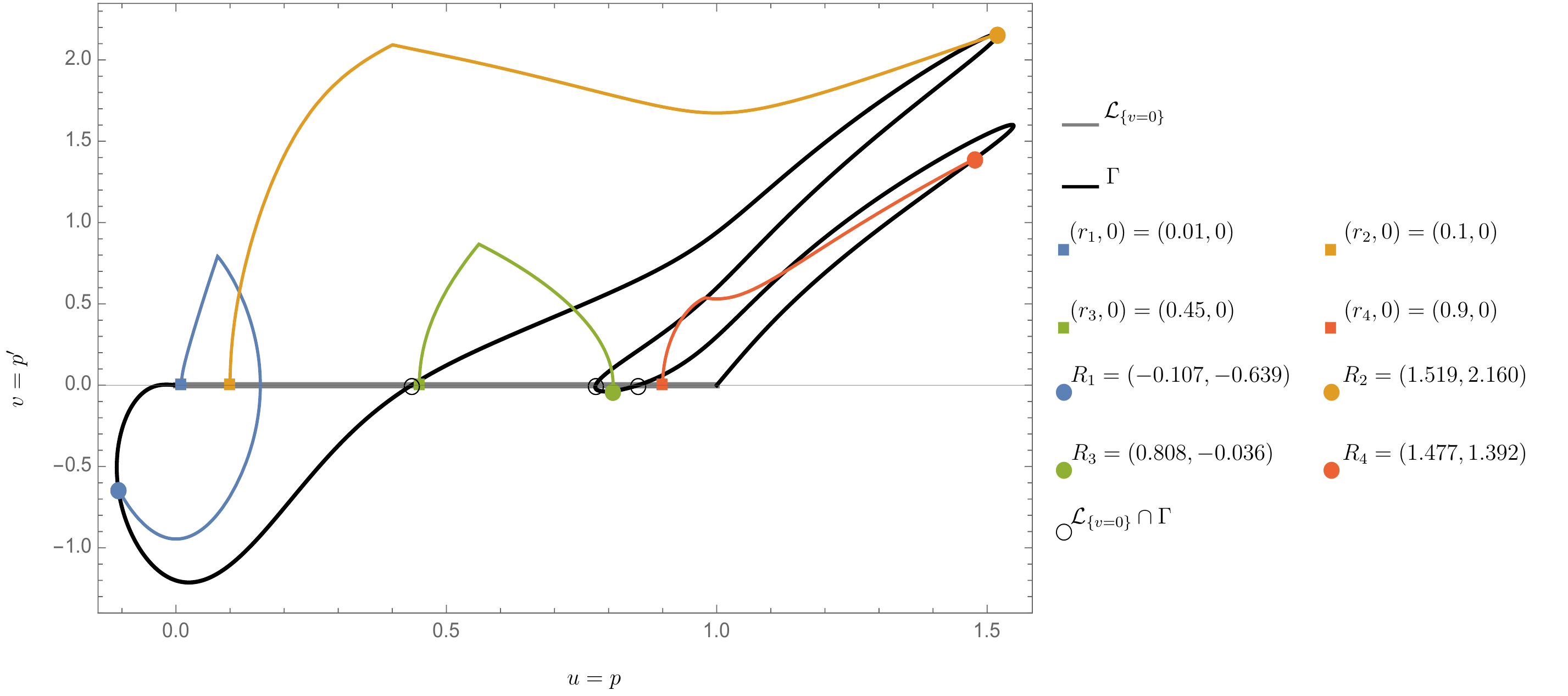}  
\caption{In the phase-plane $(u,v)$: intersection points between $\mathcal{L}_{\{v=0\}}$ and 
$\Gamma=\Phi^{\omega_{2}}_{\omega_{1}}(\mathcal{L}_{\{v=0\}})$; 
solutions of the Cauchy problem with initial conditions given by $(u(\omega_{1}),v(\omega_{1}))=(r_{i},0)$ and 
numerical approximation of the values $R_{i}=(u(\omega_{2}\, ;\omega_{1},(r_{i},0)),v(\omega_{2}\, ;\omega_{1},(r_{i},0)))$ with $i=1,\dots,4.$
The problem's setting $p''+\lambda w(x)f(p)=0$ is defined as in Proposition~\ref{prop-2}.}
\label{fig-4}
\end{figure}

In this case, the features of the nonlinear term along with the joint action of the indefinite weight 
give rise to an involved deformation of the segment $\mathcal{L}_{\{v=0\}}$. 
Nevertheless, looking at Figure~\ref{fig-4}, we can see that there exist more than one intersection point between the continuum
$\Gamma$ and the $u-$axis such that their abscissa is contained in the interval open $]0,1[.$

This way, the previous observation suggests us the following analysis. By choosing the values
 $r_{1}=0.01$, $r_{2}=0.1$, $r_{3}=0.45$ and $r_{4}=0.9$ we compute the points $R_{i}$ for $i=1,\dots,4$. 
All the results achieved are truncated at the third significant digit and so we obtain 
${R_{1}}^{v}=-0.639<0$, ${R_{2}}^{v}=2.160>0$, ${R_{3}}^{v}<-0.036$ and ${R_{4}}^{v}=1.392>0$.
The numerical details are thus represented in Figure~\ref{fig-4}.

\begin{figure}[htb]
	\centering
	\includegraphics[width=1\textwidth]{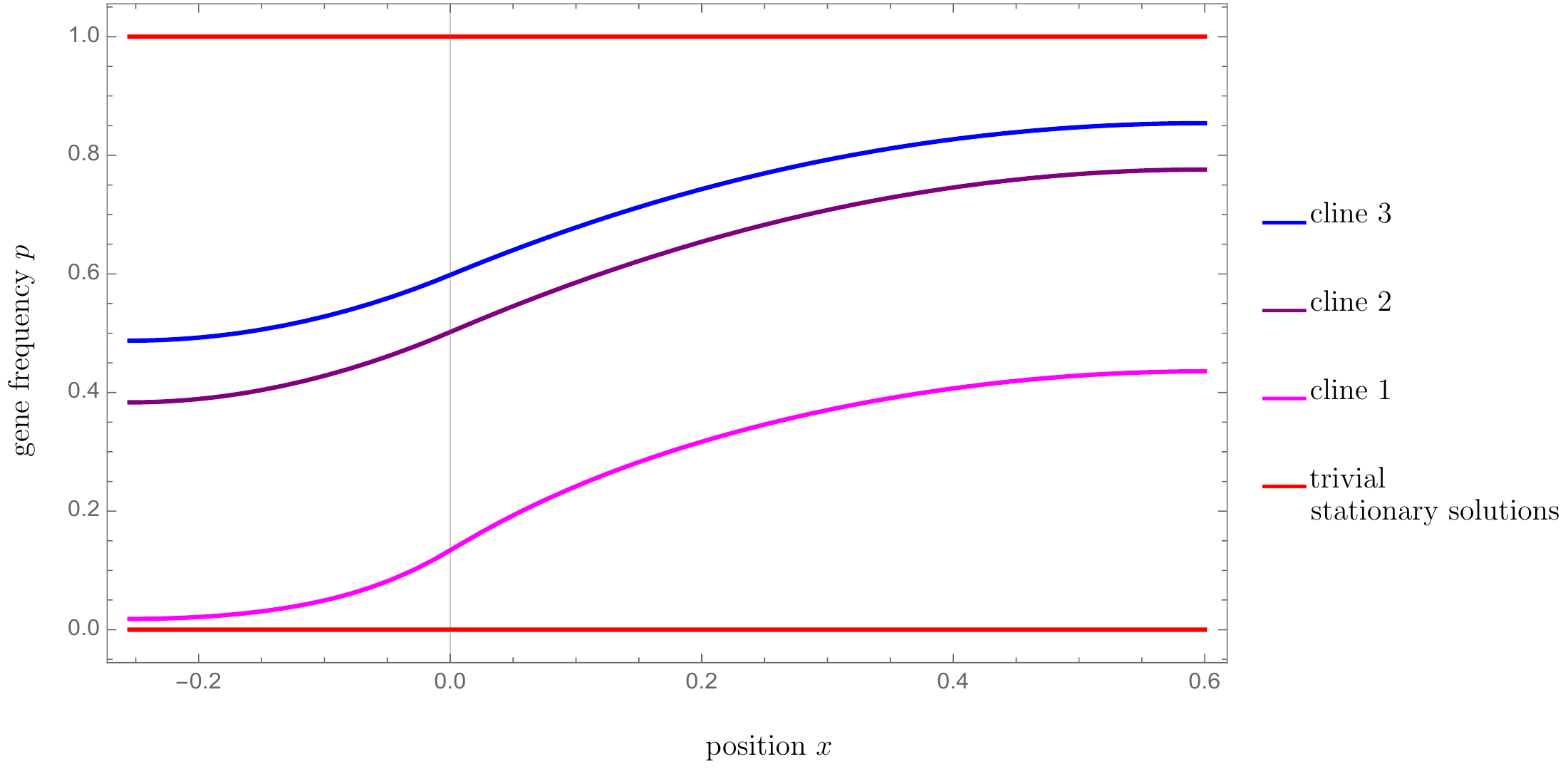}  
\caption{Polymorphic stationary solutions (clines) and trivial stationary solutions ($u\equiv0$ and $u\equiv1$) for equation \eqref{eq-rd}, 
found as positive solutions of the Neumann problem $p''+\lambda w(x)f(p)=0$ satisfying the framework of Proposition~\ref{prop-2}.}
\label{fig-5}
\end{figure}

At this point, an application of the Intermediate Value Theorem guarantees the existence of at least three initial conditions 
$(c_{j},0)$ with $j=1,2,3$, such that each respective solution of the initial value problem \eqref{eq-ivp} is also a positive solution 
of the Neumann problem \eqref{eq-3.2} we are looking for. 
Indeed, the values $c_{1}=0.436$, $c_{2}=0.776$ and $c_{3}=0.854$ satisfy the conditions in \eqref{eq-cond1}.
Finally, we display the approximation of the three nontrivial stationary solution $p(\cdot)$
of equation \eqref{eq-rd} in Figure~\ref{fig-5}. 

\section{Discussion}
\def\theequation{5.\arabic{equation}} \makeatother
\label{section-4}
When the selection gradient is described by a piecewise constant coefficient,
we have studied the conjecture proposed in \cite{LoNa-02} within a finite and one-dimensional environment. 

Summing up: we have found multiplicity of positive solutions
for two different indefinite Neumann problems defined as in \eqref{eq-3.1} where the nonlinear term is an application
$f:[0,1]\to\mathbb{R}$ which assumes two particular forms, the one in \eqref{eq-f1} or that in \eqref{eq-f2}.
In our examples, the nonlinearity $f$ is a function of class $C^2$ such that
\begin{equation*}
f(0) = f(1) =0, \quad f(s) > 0 \;\; \forall\,s\in \,]0,1[\,, \quad f'(0) > 0 > f'(1),
\leqno{(f_{*})}
\end{equation*}
and
\begin{equation*}
\emph{\text{$f$ is not concave, $s\mapsto f(s)/s$ is strictly decreasing}.}
\leqno{(H)}
\end{equation*}
Hence, uniqueness of positive solutions in general is not guaranteed for indefinite Neumann problems 
whose nonlinear term $f$ is a function satisfying $(f_{*})$ and $(H)$, and
the indefinite weight $w$ is defined on a bounded domain $\Omega$ with $\int_{\Omega}w(x)\,dx<0$.
Nevertheless, in the case of the family of functions $f_{k}(s)=s(1-s)(1+k-2ks)$
with $k\in[-1,1]$, there are issues in this direction that have not yet been answered. 

So, a question still open is the following: 
\emph{under the action of gene flow, what is the minimal set of assumptions 
under which a selection gradient will maintain a unique gene frequency cline}?

The delicate matter of the comparison between the concavity 
versus a condition about monotonicity arises also in other context than the Neumann one.
With this respect, indefinite weight problems under Dirichlet boundary conditions 
have been considered in \cite{SoZa-15} where an example of multiplicity of positive solutions was given.
As far as we know, the mathematical literature lacks of a rigorous multiplicity result in both the two cases.
This way, we see how these issues deserve to be studied in deep for both a natural mathematical 
and genetical interest.

\section*{Acknowledgments}
\noindent
I thank prof. Reinhard B\"{u}rger for inspiring me to work on this problem.  
I am also very grateful to profs. Fabio Zanolin, Carlota Rebelo and Alessandro Margheri 
for providing useful comments 
that greatly improved the manuscript.

\nocite{*}
\bibliographystyle{elsart-num-sort}
\bibliography{ReferenceSovranoElisa.bib}

\end{document}